\documentclass[oneside]{amsart}
\usepackage{amsfonts,latexsym,amssymb,amsmath,amscd}
\usepackage[english]{babel}
\usepackage[dvips]{graphicx}
\usepackage{enumerate}
\usepackage{amscd}

\newtheorem{teo}{Theorem}

\newtheorem{prop}{Proposition}

\newcommand{\be}{\begin{enumerate}}
\newcommand{\ee}{\end{enumerate}}

\newcommand{\ds}{\displaystyle}
\newcommand{\no}{\noindent}
\newcommand{\tx}{\textrm}

\newenvironment{prueba}{\no \textbf{Proof:}\:\:}{\hfill \mbox{$\square$}\par}
\newcommand{\R}{{\mathbb R}}
\newcommand{\PC}{{\mathbb P}_{\C}}

\newcommand{\N}{{\mathbb N}}
\newcommand{\Q}{{\mathbb Q}}
\newcommand{\C}{{\mathbb C}}

\newcommand{\eps}{\varepsilon}

\DeclareMathOperator{\diff}{Diff}

\newcommand{\parcial}[2]{\ds\frac{\partial #1}{\partial #2}}

\newcommand{\F}{\mathcal{F}}

\renewcommand{\H}{\mathcal{H}}

\begin{document}
\title{On the Quasi-ordinary cuspidal foliations in $(\C^3,0)$}
\author{Percy Fern\'andez-S\'anchez \and Jorge
Mozo-Fern\'{a}ndez}
\address{Percy Fern\'{a}ndez S\'{a}nchez: Instituto de Matem\'{a}tica y Ciencias Afines \\ Universidad
Nacional de Ingenier\'{\i}a y
Pontificia Universidad Cat\'{o}lica de Per\'{u} \\ Casa de las
Trece Monedas \\ Jr. Ancash 536 \\ Lima 1 \\ Per\'{u}}
\email{pefernan@pucp.edu.pe}
\address{Jorge Mozo Fern\'{a}ndez: Depto. Matem\'{a}tica Aplicada \\
University of Valladolid \\ETS de Arquitectura \\ Avda. Salamanca, s/n \\
47014 Valladolid \\ Spain} \email{jmozo@maf.uva.es}
\thanks{Both authors partially supported by CONCYTEC (Peru) under the research project 551-OAJ (2004), and by Junta de
Castilla y Le\'{o}n (Spain) under project VA123/04.}
\date{\today}
\maketitle

%\begin{abstract}
%.
%\end{abstract}

%-------------------------------------------------------------------------------------------

%  introduction and motivation

%-------------------------------------------------------------------------------------------

\section{Introduction and motivation}

We would like to study  the reduction of the singularities and the
analytic classification, in some  cases that we shall describe, of
germs of singular holomorphic foliations in $(\C^3,0)$, with
non-zero linear part. Consider, more generally, $\omega$ a germ in
$(\C^n,0)$ of an integrable $1$-form, and let

\[
\omega=\omega_1+\omega_2+\cdots
\]

\no be its decomposition in homogeneous forms
($\omega_i=\sum_{j=1}^{n}A_{ij}dx_i$, $A_{ij}$ homogeneous
polynomial of degree $i$). Suppose, moreover, that $\omega_1
 \not\equiv 0$. In general, we can write

\[
\omega_1=\sum_{i,j}^n c_{ij}x_jdx_{i}, \,\,\,\, c_{ij}\in \C.
\]

\no The integrability condition $\omega \wedge d\omega=0$ implies
that $\omega_1\wedge d\omega_1=0$. Let $C$ be the matrix
$(c_{ij})_{i,j} \in M_{n\times n}(\C)$. Writing down explicitly the
integrability condition, the coefficient of $dx_i \wedge dx_j \wedge
dx_k$ ($i<j<k$) in $\omega_1\wedge d\omega_1$ is

\[
c_i(c_{kj}-c_{jk})-c_j(c_{ki}-c_{ik})+c_k(c_{ji}-c_{ij}),
\]
where $c_i=\sum_{j=1}^n c_{ij}x_j$. Two cases appear:
\begin{enumerate}
\item $C$ is a symmetric matrix. \item $C$ is not symmetric. So,
it exists $(j,k)$ with $c_{kj}\ne c_{jk}$.
\end{enumerate}
In the last case, the polynomials $c_i$, $c_j$, $c_k$ are linearly
dependent for every $i$, $j$, $k$, and so $rk(C)\leq 2$. Moreover,
$d\omega_1(0)=d\omega(0)\ne 0$, so we are in presence of a
Kupka-type phenomenon and, in fact, it exists a biholomorphism $f$
such that $f^*\omega\wedge \eta=0$, where $\eta$ is a form in
$2$-variables. For bidimensional phenomena, lots of work have been
done.

We then focus on the symmetric case. A linear change of
coordinates changes $C$ in $P^tCP$, $P$ invertible, so we can
suppose $C$ diagonal and moreover

\[
\omega_1=\sum_{i=1}^r x_idx_i,\,\,\,\, r\leq n.
\]

If $r=n$, G. Reeb, in his thesis \cite{R} shows that there always
exists a holomorphic first integral. The behaviour of the foliation
is, then, the behaviour of a function. Using Malgrange's singular
Frobenius theorem \cite{Ma} we recover this result.

If $r<n$, some work was done by R. Moussu \cite{Mo} under additional
hypothesis. The fundamental paper of Mattei and Moussu \cite{MM}
completes the mentioned results. Let us recall in this case,
briefly, the 2-dimensional situation. The foliations studied are
defined by 1-forms $ydy+\cdots $. Following Takens \cite{T} such a
foliation has a formal normal form

\[
\omega_N=d(y^2+x^n)+x^pU(x)dy
\]

\no where $n\geq 3$, $p\geq 2$, $U(x)\in \C[[x]]$, $U(0)\ne 0$.

The generic case ($n=3$) was studied by Moussu \cite{Mo2} and a
generalization ($n\geq 3$, $2p>n$) by Cerveau and Moussu
\cite{CM}. In both cases, the reduction of the singularities of
$\omega$ (and $\omega_N$) agrees with the reduction of the curve
$y^2+x^n=0$. Projective holonomy classifies and, generically,
there is a rigidity phenomenon formal / analytic. If $n$ is even
and $2p=n$, it has been studied by Meziani \cite{Me} under some
restrictions on the values of $U(0)$. If $2p<n$
 the study was done (not in full
generality) by Berthier, Meziani and Sad \cite{BMS}. We shall call
``cuspidal'' to these foliations.

The objective of the present work is to generalize this situation to
dimension three. We want to study foliations whose linear part is
given by $d(x^2+y^2)$ or by $d(z^2)$. In this paper we shall focus
in the case $d(z^2)$. A surface that controls the resolutions of the
singularities, with an equation $z^2+\cdots=0$ will appear in the
considered cases.

Let us recall some results about reduction of the singularities of
a complex surface, following Hironaka \cite{H}. A surface $X$ in
$\C^3$ has an equation

\[
f=f_{\nu}+f_{\nu+1}+\cdots=0
\]
where $f_i$ is homogeneous of degree $i$. For such a surface,
define, at the origin:

\begin{enumerate}
\item The tangent cone, $C_{X}$, as the cone $f_{\nu}=0$.

\item The Zariski tangent cone $T_X$, as
$Spec(\mathcal{M}/\mathcal{M}^2)$, $\mathcal{M}$ being the maximal
ideal corresponding to the origin of $\C[[x,y,z]]/(f)$. This is
the smallest linear space containing $C_X$.

\item The strict tangent cone $S_X$, as the largest linear
subspace $T$ of $T_X$ such that $C_X=C_X+T$. The codimension of
$S_X$ is the minimum number of variables required to write down the
equations of $C_X$.
\end{enumerate}

The resolution of singularities of an analytic surface $X$ is a
problem that may stated as follows: to find a non-singular surface
$\widetilde{X}$ and a birational morphism $\widetilde{X}\to X$
composed of quadratic (point blow-ups) and monoidal (curve blow-ups)
transformations. These must be done in a precise order. The main
case to consider is when the three tangent spaces defined above
coincide, and the most difficult case is when, moreover, $\dim
S_X=2$. In this case, the tangent cone can be written as $z^{\nu}$.
The resolution may be controlled by Hironaka's characteristic
polyhedra of the singularities \cite{H}. The precise sequence of
blow-ups needed can be read in the polyhedra.

A kind of surface singularities whose resolution is particularly
simple, and combinatorial, are quasi-ordinary singularities. To
define them, consider a finite  projection $X \stackrel{\pi}{\to}
\C^2$ and let $\Delta$ be discriminant locus of $\pi$ (i.e. the
projection of the apparent contour). If $\Delta$ has normal
crossings the singularities of $X$ are called quasi-ordinary.

Quasi-ordinary singularities are studied not only because they are
relatively simple, but because they arise in the Jungian approach to
desingularization. First of all desingularize the discriminant locus
in order to obtain quasi-ordinary singularities. Then, the problem
(simpler) is to reduce the singularities of a quasi-ordinary
surface. Some good references of this are the articles of Giraud
\cite{G} and  Cossart \cite{Co}.

Quasi-ordinary singularities can be parametrized by fractional
power series, as branches of curves:
\[
\left\{\begin{array}{l}
x=x \\
y=y \\
z=\sum_{i,j}c_{ij}x^{\frac{i}{n}}y^{\frac{j}{n}}
\end{array}\right.
\]
By the condition of the discriminant, it can be seen that the set
of points $\{(i,j)\in \R^2:c_{ij}\neq 0\}$ is contained in a
quadrant $(a,b)+\R_{+}^2$, where $c_{ab}\neq 0$. Characteristic
pairs may be defined for this parametrization, as is the case of
curves, and they still determine the local topology of the
singularity, while the converse is not known \cite{Li}.

Coming back to foliations, this is related with the case we shall
study. More precisely, we search a class of foliations in $(\C^3,0)$
whose reduction process can be read in a quasi-ordinary surface. For
the case considered $\omega_1=d(z^2)$, by Weierstrass preparation
theorem and Tschirnhausen transformations we find that, in
appropriate coordinates, the surface is $z^2+\varphi(x,y)=0$, that
is not necessarily a separatrix. The natural generalization of
cuspidal foliations will be those with an equation
\[
\omega=d(z^2+\varphi(x,y))+A(x,y)dz.
\]

In fact, in a recent work, Frank Loray \cite{L}  finds an analytic
normal form as

\[
\omega=dF+zdG+zdz,
\]

\no where $F,G\in \C\{x,y\}$, for integrable holomorphic foliations
with linear part not tangent to the radial vector field. Note that a
coordinate change $z\to z-G(x,y)$ in Loray's form gives an equation
like our expression for the foliations. This is integrable if and
only if $d\varphi \wedge dA=0$, i.e. if $\varphi$, $A$ are
analytically dependent. As are shall restrict to the quasi-ordinary
case, we have that $\varphi(x,y)=x^py^qU(x,y)$, with $U$ a unit. A
convenient change of variable in $x$, $y$, allows us to suppose that
$\varphi(x,y)=x^py^q$. Let $d=\gcd (p,q)$ $p=dp'$ $q=dq'$. The
integrability condition $d\varphi\wedge dA=0$ is then that
$A(x,y)=L(x^{p'}y^{q'})$ where $L(u)\in \C\{u\}$.

The plan of this paper is as follows. In section 2 we shall review
the notion of simple singularity of a foliation, in the sense
defined by Cano and Cerveau \cite{CC}, and its analytic
classification according to Cerveau and Mozo \cite{CeMo}. Section 3
is devoted to describe the resolution of singularities of the
quasi-ordinary foliations we are going to study, and the topology of
the exceptional divisor. In section 4, we construct a Hopf fibration
associated to the quasi-ordinary foliations, making a reduction of
the separatrix to a canonical form. Finally, section 5 is devoted to
present the main result of the paper: In the considered cases, the
holonomy of a certain component of the exceptional divisor
classifies analytically the foliation. The cases we study, as we
shall see, are essentially the sames that are studied in dimension
two.

Some notations used throughout the paper are presented here.
$\diff (\C ,0)$ will denote the group (under composition) of germs
of analytic diffeomorphisms of $(\C,0)$. If $\Omega$ denotes a
holomorphic integrable 1-form, defining a foliation, and $D$ is a
component of the divisor obtained after reduction of
singularities, ${\mathcal H}_{\Omega, D} : \pi_{1}(D\setminus
{\mathcal S})\longrightarrow \diff (\C,0)$ is the holonomy
representation, defined over a transversal to $D$ (omitted from
the notation), where ${\mathcal S}$ is the singular set of the
reduced foliation.

{\bf Acknowledgements.-} This paper was begun during the stay of the
first author at the University of Valladolid and was finished at the
Instituto de Matem\'{a}tica y Ciencias Afines (IMCA), of Lima, and
the Pontificia Universidad Cat\'{o}lica de Per\'{u}. We want to
thank these institutions.

%------------------------------------------------------------------------------------

% Singularidades Simples

%-------------------------------------------------------------------------------------------

\section{Simple singularities of foliations and analytic classification}

The process of reduction of singularities for a holomorphic
foliation is well known in dimension two. After a finite number of
point blow-ups performed in any order, a germ of analytic space
and a foliation are obtained, and around the singular points, the
foliation is generated by a one-form
\[
\omega=(\lambda y+h.o.t.)dx + (\mu x+h.o.t.)dy,
\]
with $\mu\neq 0$, $\frac{\lambda}{\mu}\not \in \Q_{<0}$.

The analytic classification is well studied in a wide variety of
cases:

\begin{enumerate}
\item If $\frac{\lambda}{\mu}\not \in \R_{\geq 0}$, $\omega$ is
analytically linearizable, i. e., there exists and analytic
diffeomorphism $\phi:(\C^2,0)\to (\C^2,0)$ such that
\[
\phi^*\omega \wedge (\lambda ydx+\mu xdy)=0
\]

\item If $\frac{\lambda}{\mu}\in \R_{>0} \setminus \Q$, but it is
not ``well-approached" by rational numbers, it is also
linearizable. If it is well-approached, we face a problem of small
divisors, and the situation becomes more complicated.

\item If $\lambda=0$ or $\frac{\lambda}{\mu}\in \Q_{+}$, Martinet
and Ramis find a large moduli space formal/analytic. In this case
the classification of the foliation agrees with the classification
of the holonomy of a strong separatrix (i.e. a separatrix in the
direction of a non-zero eigenvalue). Moreover in the resonant case
($\frac{\lambda}{\mu}\in \Q_{+}$) or in the saddle-node case
($\lambda=0$) with analytic center manifold, the conjugation of
the foliation is fibered. This means the following: choose
coordinates $x$, $y$ such that  the axis are the separatrices,
$y=0$ being a strong one; the foliations are defined by 1-forms
\[
\omega_i=yA_i(x,y)dx+\mu /\lambda x(1+B_i(x,y))dy,
\]
with $i=1,2$. Let $h^{(i)} (x)$ be the holonomies of $y=0$, supposed
conjugated. Then the foliations are conjugated by a diffeomorphism
$\phi(x,y)=(x,yg(x,y))$.
\end{enumerate}

The singularities obtained after this reduction process are called
simple or reduced. The class of simple singularities is stable
under blow-ups. Let us observe that the notion of simple
singularity is not only analytic, but formal: if $\omega_1$,
$\omega_2$ are analytic $1$-forms, and $\hat{\phi}$ is a local
diffeomorphism such that $\hat{\phi}^{*}\omega_1\wedge
\omega_2=0$, then $\omega_1$ has a simple singularity if and only
if $\omega_2$ has.

If the dimension of the ambient space is greater or equal than
three, the notion of simple singularity has been developed in
\cite{CC}, \cite{Ca}, and its analytic classification studied in
\cite{CeMo}. The reduction of singularities is only achieved when
the dimension of the ambient space is at most three, and in this
case, simple singularities are the final ones obtained after the
reduction process. Let us summarize here, for convenience of the
reader, the main results in dimension three.

First of all, let us recall the notion of ``dimensional type". A
foliation has dimensional type $r$ if there exist analytic (resp.
formal) coordinates such that the foliation is defined by an
integrable $1$-form $\omega$ that can be written in coordinates
$x_1,\cdots,x_r$ ($r\leq n$), but not less. So, a
three-dimensional singularity of foliation has dimensional type
$2$ and $3$. For instance, if we are in presence of a Kupka
phenomenon, the dimensional type is $2$. The notions of formal
dimensional type or analytic dimensional type are equivalent, as
seen  in \cite{CeMo}. So, we have simple singularities of
dimensional types $2$ and $3$. If the dimensional type is $2$,
simple singularities are defined by a simple $2$-dimensional
$1$-form. They have $2$ separatrices, of which at most one is
formal.

If the dimensional type is tree, simple singularities are the ones
that admit one of the following formal normal forms:
\begin{align}
\omega=xyz\left(\alpha \frac{dx}{x}+\beta \frac{dy}{y} +\gamma
\frac{dz}{z}\right),
\end{align}
with $\frac{\alpha}{\gamma},\:\frac{\beta}{\gamma},\:
\frac{\alpha}{\beta}\not \in \Q_{-}$ (and $\alpha\beta\gamma\neq 0$,
as the dimensional type is $3$). This is the linearizable case. If,
for instance, some of the quotients is not real, the linearization
is analytic \cite{CL}.

\begin{align}
\omega_N=xyz\left( x^py^qz^r \right)^s\left[\alpha
\frac{dx}{x}+\beta \frac{dy}{y}+
\left(\lambda+\frac{1}{(x^py^qz^r)^s}\right)
\left(p\frac{dx}{x}+q\frac{dy}{y}+ r\frac{dz}{z}\right)\right],
\end{align}
where $p,q,r\in \N$, $qr\neq 0$, $s\in \N^*$, $\alpha$, $\beta$
constants, not both zero. This is the \emph{resonant case}.
Several things can be said about foliations that are formally
equivalent to this normal form:

\begin{enumerate}
\item $\F$ has three separatrices, of which at most one is formal
(which, in the preceding coordinates,  would be $x=0$). This is a
confluence of simple two-dimensional singularities defined along the
axis. Saddle-nodes only appear if $p=0$, and only in this case  the
existence of a formal, non convergent separatrix is possible.

\item The holonomy group of $z=0$ (strong separatrix) classifies
analytically the foliation. Moreover, the conjugations is fibered
if the three separatrices are convergent.

\item If $\frac{\alpha}{\beta}\not \in \Q$, there is a rigidity
phenomenon: every such foliation is analytically equivalent to
$\omega_N$.
\end{enumerate}

A typical case in which we are in presence of a simple singularity
and that will appear in the sequel, is when the foliation is defined
by a 1-form
\begin{align}\label{simple}
\omega=xyz\left[(p+A(x,y,z))\frac{dx}{x}+(q+B(x,y,z))\frac{dy}{y}+
(r+C(x,y,z))\frac{dz}{z}\right],
\end{align}
with $p,q,r\in \N^*$,
$\nu(A), \nu(B), \nu(C)>0$.

More can be said: the transformation $\phi$ that converts $\omega$
in its formal normal form $\omega_N$, even if it is not analytic, it
is transversally formal and fibered. This means in particular that
such a  $\phi$ can be found in the form
\[
\phi(x,y,z)=(x,y,\varphi(x,y,z)).
\]

The existence of local holomorphic first integrals, according to
Mattei and Moussu \cite{MM}, is equivalent to the periodicity of
the holonomy group. Moreover, an integrable 1-form $\omega$, that
generates a reduced foliation of dimensional type three, has a
holomorphic first integral if and only if there exists analytic
coordinates $(x,y,z)$ such that
$$
\omega \wedge (pyz dx+ qxz dy+r xy dz) =0,
$$
where $p,q,r\in \N^\ast$.

%------------------------------------------------------------------------------------

% Desingularizacion

%-------------------------------------------------------------------------------------------
\section{Reduction of singularities and topology of the
 divisor}\label{section3}

In this paper, we shall study the analytic classification of
quasi-ordinary cuspidal foliations in dimension three, i.e.,
foliations such that, in appropriate coordinates, can be defined
by an integrable 1-form
$$
\omega= d(z^2+x^py^q)+A(x,y)dz.
$$
The integrability condition here is equivalent to $d(x^py^q)\wedge
dA=0$. So, let $d= gcd (p,q)$, $p=dp'$, $q=dq'$. Such a 1-form can
be written as
$$
\omega= d(z^2+x^py^q)+(x^{p'}y^{q'})^k h(x^{p'}y^{q'})dz,
$$
where $h(u)\in \C \{ u\}$, $h(0)\neq 0$. Fixing $p$, $q$, we shall call $\Sigma_{pq}$
the set of holomorphic foliations that are analitically equivalent to the foliation defined by one of these
1-forms.

As it will become clear from the development of the paper, the
separatrices of this foliation have the equation
$$
z^2+x^py^q+ h.o.t. =0,
$$
and Weierstrass preparation theorem and Tschirnhausen transformation
show that this separatrix is analytically equivalent to
$z^2+x^py^q=0$.

The reduction of singularities for these foliations is quite
simple, similar to plane curves, and it is the main objective of
this section their detailed analysis. For convenience, we divide
the problem in three cases:

\be[{Case }1.]

\item $p$, $q$ even.

\item $p$ even, $q$ odd.

\item $p$, $q$ odd.

\ee

\be[{{\bf Case }}1.] \item Suppose $p$, $q$ are  even, and
$d=2d'$. If $k>d'$, the reduction of the singularities is obtained
after $\dfrac{p+q}{2}$ blow-ups: \be \item First of all, blow up
$\dfrac{p}{2}$ times the $y$-axis. We obtain a sequence of
divisors $D_1,\ldots ,D_{p/2}$, topologically germs $(\PC^1 \times
\C , \PC^1 )$. The intersection of two consecutive components is a germ of a
line $(\C, 0)$, $L_i=D_i\cap D_{i+1}$, $1\leq i <\dfrac{p}{2}$. In
the appropriate chart, these blow-ups have the equations
$$
\left\{ \begin{array}{rcl}
x & = & x \\
y & = & y \\
t_{i-1} & = &x\cdot t_i,
\end{array}
\right.
$$
where $t_0=z$, $1\leq i \leq \dfrac{p}{2}$. \item Then blow-up
$\dfrac{q}{2}$ times the $x$-axis, obtaining again a sequence of
divisors $D_{\frac{p}{2}+1},\ldots ,D_{\frac{p+q}{2}}$,
topologically equal to $(\PC^1 \times \C , \PC^1 )$. Again, the
intersection between two consecutive components is a line $L_i=
D_i\cap D_{i+1}$, $\dfrac{p}{2}+1 \leq i < \dfrac{p+q}{2}$. Now,
the coordinates of the blow-ups are
$$
\left\{ \begin{array}{rcl}
x & = & x \\
y & = & y \\
t_{i-1} & = &y\cdot t_i,
\end{array}
\right.
$$
with $\dfrac{p}{2}<i <\dfrac{p+q}{2}$.

The result of the composition of all the blow-ups in the preceding
charts is the map $\pi (x,y,t_{\frac{p+q}{2}})=
(x,y,x^{\frac{p}{2}}\cdot y^{\frac{q}{2}}\cdot
t_{\frac{p+q}{2}})$. The pull-back of the foliation is given by
\begin{eqnarray*}
\lefteqn{\pi^\ast \omega  =  x^{p-1}y^{q-1} \cdot \left[ 2xyt dt
+(t^2+1)xy \left( p\frac{dx}{x}+q\dfrac{dy}{y}\right) + \right. }
& &
\\
 & + & \left. (x^{p'}y^{q'})^{k-d'} h(x^{p'} y^{q'})xyt\cdot
\left( \dfrac{p}{2} \dfrac{dx}{x} +\dfrac{q}{2} \dfrac{dy}{y} +
\dfrac{dt}{t} \right) \right]
\end{eqnarray*}
(here $t= t_{\frac{p+q}{2}}$).

The foliation, now, is reduced. Let ${\mathcal S}$ be the singular
locus of this reduced foliation. ${\mathcal S}$ is an analytic,
normal crossing space of dimension one, composed by:
\be
\item The
lines $L_i$ of intersection of the divisors. These are resonnant
singular points of dimensional type two.
\item The lines $L$, $L'$
in $D_{\frac{p+q}{2}}$ of equations $(y=0, t=i)$, $(y=0, t=-i)$,
and also the lines $M'$, $M''$ in $D_{\frac{p}{2}}$ of equations
$(x=0, t=i)$, $(x=0, t=-i)$ (in the last chart). These lines are
the intersections of the two separatrices $S'$, $S''$ with the
divisors.
\item The intersection $P_i:=D_{\frac{p}{2}}\cap D_i$,
$\dfrac{p}{2}<i\leq \dfrac{p+q}{2}$ is a projective line composed of
points of dimensional type two, except at the corners: \be
\item
$m_i=P_i\cap L_i= D_{\frac{p}{2}}\cap D_i \cap D_{i+1}$,
$\dfrac{p}{2}<i<\dfrac{p+q}{2}$. These are the resonnant singular
points of dimensional type three, having $D_{\frac{p}{2}}$, $D_i$,
$D_{i+1}$ as separatrices.
\item
\begin{eqnarray*}
m' & := & D_{\frac{p}{2}}\cap D_{\frac{p+q}{2}} \cap S' = L' \cap
M' \cap P_{\frac{p+q}{2}}, \text{ and } \\
m'' & := & D_{\frac{q}{2}}\cap D_{\frac{p+q}{2}}\cap S'' = L''
\cap M''\cap P_{\frac{p+q}{2}}.
\end{eqnarray*}

These are the resonant singular points of dimensional type three
corresponding to the separatrices of the foliations.
\ee
\ee
\ee
According to the preceding description of the resolution of the
singularities, we have all the information about the topology of
$D_i\setminus {\mathcal S}$, and more precisely about the
fundamental group of these components. We have:
\begin{itemize}
\item $D_1\setminus {\mathcal S}$ is topologically $\C \times \C$,
so simply connected.
\item $D_i\setminus {\mathcal S}$
($1<i<\dfrac{p}{2}$) is topologically $\C^\ast \times  \C$. The
generator of the fundamental group is a loop $\gamma_i$ that turns
around $L_i$ (or $\gamma_i^{-1}$ around $L_{i-1}$).
\item
$D_{\frac{p}{2}+1}\setminus {\mathcal S} \cong \C^\ast \times \C$.
The fundamental group is generated by a loop $\alpha_i$ around
$P_{\frac{p}{2}+1}$.
\item $D_i\setminus {\mathcal S} \cong
\C^{\ast }\times  \C^\ast$ ($\dfrac{p}{2}+1<i<\dfrac{p+q}{2}$). The
fundamental group has generators $\gamma_i$ around $L_i$ and
$\alpha_i$ around $P_i$, that commute.
\item
$D_{\frac{p+q}{2}}\setminus {\mathcal S} \cong (\C\setminus \{
m',m''\}) \times \C^\ast$. We have one loop $\alpha_{\frac{p+q}{2}}$
around $P_{\frac{p+q}{2}}$ and loops $\gamma'$, $\gamma''$ around
the separatrices (i.e., around $m'$, $m''$). \item
$D_{\frac{p}{2}}\setminus {\mathcal S} \cong \C^2\setminus {\mathcal
C}$, where ${\mathcal C}$ is the curve with coordinates
$t^2_{\frac{p}{2}}+y^q=0$, composed of two smooth branches that meet
tangentially at the origin. In this case (see \cite{La}), $\pi_{1}(
\C^2\setminus {\mathcal C})$ is the group, written in terms of
generators and relations as
$$
\pi_{1}(\C^2\setminus {\mathcal C})= \langle \alpha, \beta;
\alpha^{\frac{q}{2}}\beta =\beta\alpha^{\frac{q}{2}} \rangle.
$$

These loops go as follows. Consider the curve $t^2_{\frac{p}{2}}+ y^q=0$
on $\C^2$, and cut by $y=1$. You obtain
$\C\setminus \{ m', m'' \}$; then $\alpha$ is a loop in $y=1$ that
turns around these two points $m'$, $m''$, and $\beta$ is a loop
in $t_{\frac{p}{2}}=0$ that turns around the origin. At the end of
the reduction process, $\alpha$ is going to be a loop in
$D_{\frac{p}{2}}$ around the two separatrices, and $\beta$ a loop
around $P_{\frac{p+q}{2}}$ ``between $S'$ and $S''$".
\end{itemize}

The case $k=d'$ ($2k=d$) is almost identical, except for some
values of the coefficient $h(0)$. More precisely, after
$\dfrac{p+q}{2} $ blow-ups, in order to obtain the complete
reduction of singularities (i.e., simple singular points) it is
necessary and sufficient that
$$
h(0)^2\neq \frac{(16+r)^2}{16+2r},\ \forall\, r\in \Q_{>0}.
$$
Moreover, if in the preceding expression we put $r=0$, we have
then $h(0)=\pm 4$. In this case, only one separatrix is obtained,
but it is a three-dimensional saddle-node, the divisor being the
weak separatrix (then convergent). We shall assume that this is
not the case, i.e., if $k=d'$ we shall assume that
$$
h(0)^2 \neq \frac{(16+r)^2}{16+2r},\ \forall\, r\in \Q_{\geq 0}.
$$
The reader may verify that this condition is equivalent to
${\mathcal P}_2$ property in \cite{Me,Me-tesis} (i.e. $h(0)\neq
\pm 2 \left( \sqrt{r}+ \dfrac{1}{\sqrt{r}} \right) $, $\forall\,
r\in (0,1] \cap \Q$).

Suppose now that $k<d'$. In this case, the reduction of
singularities is achieved blowing-up $kp'$ times the $y$-axis and
$kq'$ times the $x$-axis. After these, in the last chart we obtain
as singularities the sets $L'= (x=t=0)$, $M'=(y=t=0)$, $L''= (x=0,
t=1)$, $M''= (y=0, t=1)$. These are also two singular points of
dimensional type three, namely $m'= L'\cap M'\cap P_{k(p'+q')}$,
$m''=L''\cap M'' \cap P_{k(p'+q')}$ (with analogous notations as
before), corresponding respectively to the points $(0,0,0)$ and
$(0,0,1)$. But now $m''$ is a saddle-node, so the separatrix $S''$
is maybe formal. In this paper, we shall assume that always $S''$
is convergent, i.e., there is a center manifold.

\item Suppose $p$ even, $q$ odd. If $k>d$, the reduction of
singularities is obtained after the following sequence of
blow-ups. \be \item First, blow-up $p/2$ times the $y$-axis,
obtaining divisors $D_1,\ldots ,D_{\frac{p}{2}}$ linked by lines
$L_1,\ldots ,L_{\frac{p}{2}-1}$. The equations of these blow-ups
are
$$
\left\{
\begin{array}{rcl}
x & = & x \\
y & = & y \\
t_i & = & x\cdot t_{i+1},
\end{array}
\right.
$$
where $t_0:=z$, $i<\dfrac{p}{2}$.

\item Blow-up $\dfrac{q-1}{2}$ times the $x$-axis, obtaining
$D_{\frac{p}{2}+1},\ldots ,D_{\frac{p+q-1}{2}}$ joined by lines
$L_i=D_i\cap D_{i+1}$, and $D_i$ joined to $D_{p/2}$ by a
projective $P_i$. The equations are
$$
\left\{
\begin{array}{rcl}
x & = & x \\
y & = & y \\
t_i & = & y\cdot t_{i+1},
\end{array}
\right.
$$
$\dfrac{p}{2}\leq i <\dfrac{p+q-1}{2}.$

\item It appears a tangency in the singular locus. In order to
break it, blow-up again the $x$-axis and take a chart centered in
the point corresponding to $t_{\frac{p+q-1}{2}}$. The equations
are now
$$
\left\{
\begin{array}{rcl}
x & = & x \\
y & = & s\cdot t_{\frac{p+q-1}{2}} \\
t_{\frac{p+q-1}{2}} & = & t_{\frac{p+q-1}{2}},
\end{array}
\right.
$$
and we obtain a new component $D'$ such that $D'\cap
D_{\frac{p+q-1}{2}}=L_{\frac{p+q-1}{2}}$, $D'\cap
D_{\frac{p}{2}}=P'$.
\item Finally, blow-up again the $x$-axis, in
order to obtain normal crossings. We obtain a final component $D''$
and the only separatrix $S$ of the foliation cuts $D'$
transversely in a line $L$ (and $D_{p/2}$ in a line $M$). We have
$L'=D'\cap D''$ and $P''=D''\cap D_{p/2}$. \ee

The singular points of dimensional type three are
$m_i:=D_{p/2}\cap D_i\cap D_{i+1}$
($\dfrac{p}{2}<i<\dfrac{p+q-1}{2}$), $m_{\frac{p+q-1}{2}}:=
D_{p/2}\cap D_{\frac{p+q-1}{2}}\cap D'$, $m'=D_{p/2} \cap D'\cap
D''$. and $m=D_{p/2}\cap D' \cap S$.

The topology of the components is as in Case 1. If ${\mathcal S}$
is the singular locus, $D_1\setminus {\mathcal S} \cong \C^2$ is
simply connected, $D_i\setminus {\mathcal S} \cong \C^\ast \times
\C$ if $1<i<\dfrac{p}{2}$, $D_{\frac{p}{2}+1}\setminus {\mathcal
S} \cong \C^\ast \times \C$, $D_i\setminus {\mathcal S} \cong
\C^\ast \times \C^\ast$ if $\dfrac{p}{2}+1<i<\dfrac{p+q-1}{2}$,
$D'\setminus {\mathcal S} \cong (\C\setminus \{ m,m'' \} ) \times
\C^\ast$, $D''\setminus {\mathcal S}\cong \C \times \C^\ast$.
Finally, $D_{p/2}\setminus {\mathcal S}\cong \C^2 \setminus
{\mathcal C}$, where ${\mathcal C}$ is the curve with coordinates
$t_{p/2}^2 +y^q=0$. As before,
$$
\pi_1 (\C^2 \setminus {\mathcal C}) = \langle \alpha, \beta ;
\alpha^q=\beta^2 \rangle .
$$
\vspace{0.5cm}
\begin{figure}[h]
\begin{center}
\includegraphics[height=6.5cm,width=12cm]{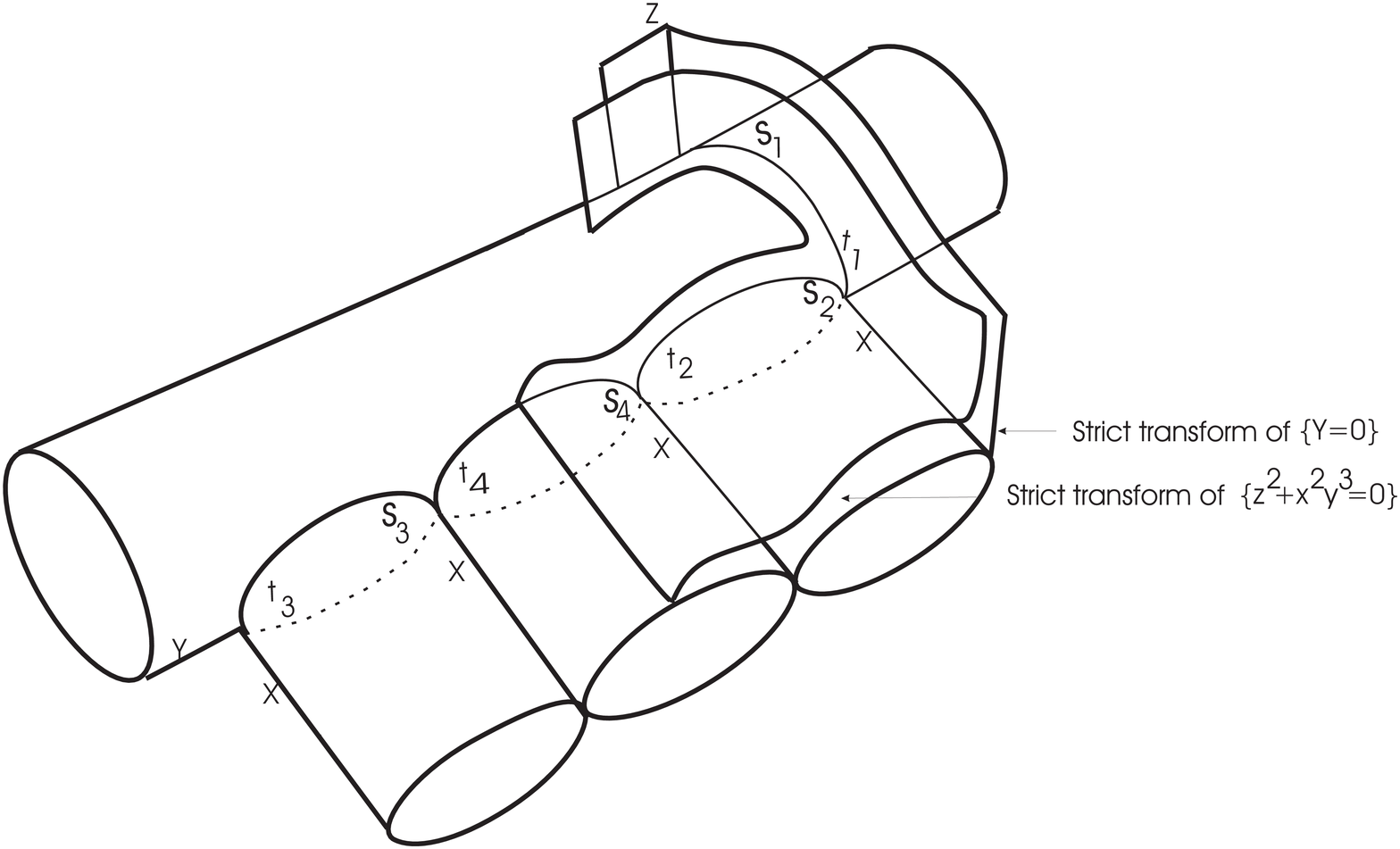}
\caption{\footnotesize The reduction of the surface
$z^2+x^2y^3=0$}. \label{reduction}
\end{center}
\end{figure}
When $k<d$, as in dimension two, the situation is as in Case 1, with $k<d'$.

\item $p$, $q$ odd. Now, the resolution is something different than before.
First, blow-up $\dfrac{p-1}{2}$ times the $y$-axis and
$\dfrac{q-1}{2}$ times the $x$-axis obtaining $D_1,\ldots
,D_{\frac{p-1}{2}}, D_{\frac{p+1}{2}},\ldots
,D_{\frac{p+q}{2}-1}$. In the new coordinates $(x,y,
t:=t_{\frac{p+q}{2}-1})$ the singular locus is given by the three
coordinate axis, that corresponds to the intersection of the
divisors and the intersection of the cone $t^2+xy=0$ with the
divisors.

Now, blow-up the origin, obtaining $P$, a projective $\PC^2$. The
three coordinate axis, now transverse to $P$, continue being
singular. Over $P$, the singular locus is composed by two
projective lines and a conic tangent to both lines. In order to
finish, blow-up twice each of the axis $x$ and  $y$ transverse to
$P$, obtaining $D_{(1)}'$, $D_{(1)}''$, $D_{(2)}'$, $D_{(2)}''$.

With respect to the topology of the divisors, the only interesting
case (i.e., not similar to the preceding ones) to comment is
$P\setminus {\mathcal S}$. As we said before, $P\cap {\mathcal S}$
is composed by two lines and a regular conic, so
$$
\pi_1 (P\setminus {\mathcal S}) \cong \langle \alpha, \beta;
\alpha^2\beta =\beta \alpha^2 \rangle.
$$
\ee
\section{Reduction of the separatrix to a canonical form}

Let $\F$ be a germ of a singular foliation defined on $(\C^3,0)$,
and let $\pi:(M,D)\to (\C^3,0)$ be the minimal reduction of the
singularities of $\F$ in Cano-Cerveau sense, as described above
\cite{CC}. Let $\tilde{\F}$ be the  strict transform of the
foliation $\F$ by $\pi$ and let $D_i$ be a component of the
exceptional divisor $D$.

 We recall, that a \emph{Hopf fibration} $\H_{\F_{\Omega}}$ adapted to
 $\F_{\Omega}\in \sum_{pq}$ is a  holomorphic transversal fibration
$f:M\to D_i$ to the foliation $\F_{\Omega}$, i.e:
\begin{enumerate}
\item $f$ is a retraction, more precisely, $f$ is a submersion and
$f|_{D_i}=Id_{D_i}$.

\item The fibers $f^{-1}(p)$ of $\H_{\F_{\Omega}}$ are contained in
the separatrices of $\F_{\Omega}$, for all $p\in D_i\cap Sing
(\tilde{\F}_{\Omega}) $.

\item The fibers $f^{-1}(p)$ of $\H_{\F_{\Omega}}$ are transversal
to the foliation $\F_{\Omega}$, for all $p\in D_i\setminus Sing
(\tilde{\F}_{\Omega})$.
\end{enumerate}

We shall be interested in finding a Hopf fibration adapted to the
foliation, relative to a particular component of the exceptional
divisor. For, if $p$ is even, call $\tilde{D}:= D_{p/2}$, i.e.,
the last component obtained after the first sequence of line
blow-ups. If $p$ and $q$ are odd, $\tilde{D}:= P$, i.e., the
projective obtained after the (only) point blow-up.

The     task of finding a Hopf fibration associated to the
foliation ${\F}_{\Omega}$ is not easy in the actual coordinates
$(x,y,z)$. As it is done in the two-dimensional case, to overcome
this obstacle, we analyze the desingularization of
${\F}_{\Omega}$
in order to obtain a simple equation for the separatrices.

From  Section (\ref{section3}) we know that the foliation
$\F_{\Omega}\in \Sigma_{pq}$ defined by the one-form
\[
\Omega=d(z^2+(x^{p'}y^{q'})^d)+(x^{p'}y^{q'})^kh(x^{p'}y^{q'})dz,
\]
has a separatrix analytically equivalent to $S:
z^2+(x^{p'}y^{q'})^r=0$ for some $r\in \N$. In order to find a
Hopf fibration $\H_{\F}$ of the foliation $\F$, we need to
normalize the one-form $\Omega$ such that the foliation defined by
this normal form  has exactly  $S: z^2+(x^{p'}y^{q'})^r=0$ as
separatrix, for certain $r$. So, the strict transformed of $S$ by
the desingularization is an hyperplane in these coordinates and
invariant by Hopf fibration.

\begin{prop}\label{normal} The foliation $\F_{\Omega}$ is analytically
equivalent to a foliation defined by the one-form
\[
d(z^2+(x^{p'}y^{q'})^r)+g(x^{p'}y^{q'},z).x^{p'}y^{q'}z\left(2\frac{dz}{z}-p'\frac{dx}{x}-q'\frac{dy}{y}\right),
\]
where  $r=d$ if $2k\geq d$ and $r=2k$ if $2k< d$. In particular,
the separatrix of the foliation $\F_{\Omega}$ is analytically
equivalent to $S:\:z^2+(x^{p'}y^{q'})^r=0$.
\end{prop}

\begin{prueba}  The foliation $\F_{\Omega}$ is defined by the
$1$-form
\[
\Omega=d(z^2+(x^{p'}y^{q'})^d)+(x^{p'}y^{q'})^kh(x^{p'}y^{q'})dz,
\]
where $(p,q)=d$, $p=p'd$, $q=q'd$. That is,
$\Omega$ is the pull-back of the 1-form
$\Omega_0=d(z^2+u^d)+u^kh(u)dz$ by the ramified fibration
\[
\begin{array}{ccl}
\rho: (\C^3,0)&\to& (\C^2,0) \\
(x,y,z) &\to & (x^{p'}y^{q'},z)=(u,z).
\end{array}
\]

The equation of the separatrices of $\Omega$ is of the form $z^2+u^d+h.o.t. =0$ (if $d$ is even, this is a joint equation, i.e.,
the product of the two separatrices).

Using Weierstrass' preparation theorem, we can assume that the local
equation of the separatrix is a polynomial in $z$:
$z^2+a(u)z+b(u)=0$, with $a(0)=b(0)=0$. If $\Phi_1
(u,z)=(u,z-\frac{a(u)}{2})$ is the Tschirnhausen transformation,
then the pull-back $\Phi_1^\ast \Omega_0$ has $z^2+c(u)=0$ as
separatrix, with $c(u)=b(u)-\frac{a(u)^2}{4}=u^r f(u)$, $f(0)\neq
0$. If $d>2$ (cuspidal case), we have that $\nu (a)>1$, $\nu (b)>2$,
and then $r>2$. In fact, $r=d$ when $2k>d$ or $r=2k$ when $2k\leq d$
(see \cite{CM,BMS,Ca}). Similar computations are valid when $d=1$ or
$d=2$ (in these cases, $2k\geq d$).

Let us write this reduced equation of the separatrices as
$$
\frac{z^2}{f(u)}+u^r=0,
$$
and let $f(u)^{1/2}$ be a square root of the unit $f(u)$. If $ \Phi_2 (u,z)= (u,z\cdot f(u)^{1/2})$, and $\Phi:= \Phi_1 \circ \Phi_2$, then $\Phi^\ast \Omega_0$ has $z^2+u^r=0$ as separatrix. This map has the form
$$
\Phi (u,z)= \left(
u,z\cdot f(u)^{1/2}-\frac{a(u)}{2}
\right) .
$$

Consider the diagram
$$
\begin{CD}
\C^3 @>{\rho}>> \C^2 \\
@VFVV @VV{\Phi}V \\
\C^3 @>{\rho}>> \C^2
\end{CD}.
$$

We want to find a diffeomorphism $F= (F_1,F_2,F_3)$ that makes commutative the diagram, i.e., that

\[
(F_1^{p'}F_2^{q'},F_3)=\left(x^{p'}y^{q'},zf(x^{p'}y^{q'})^{\frac{1}{2}}-\frac{a(x^{p'}y^{q'})}{2}\right).
\]
For, we may choose $F_1=x$, $F_2=y$, $F_3=z\cdot f(x^{p'}y^{q'})^{1/2}- \dfrac{a(x^{p'}y^{q'})}{2}$.
The form $\Phi^* \Omega_0$, having $z^2+u^r=0$ as a separatrix is,
up to a unit, $d(z^2+u^r)+g(u,z)(2udz-dzdu)$, so $F^*\Omega_0$
defines the same  foliation that
\[
d(z^2+(x^{p'}y^{q'})^r)+g(x^{p'}y^{q'},z).x^{p'}y^{q'}z\left(2\frac{dz}{z}-p'\frac{dx}{x}-q'\frac{dy}{y}\right).
\]
We reproduce part of the proof presented in \cite{CM} in order to
find the transformation $F$ fibered.
\end{prueba}

As a  consequence of this normal form for $\F_{\Omega}$, there
exists coordinates $(x,y,z)$, such that the separatrix $S$ of the
normal form is given by the equation: $z^2+(x^{p'}y^{q'})^r=0$,
where $r$ is as in the Proposition (\ref{normal}), and not only
``analytically equivalent to". Now, we can find a Hopf fibration,
from a holomorphic vector field $X_1$ for which $S$ is an invariant
set, that is
\[
X_1=\left\{\begin{array}{ll}
x\parcial{}{x}+\dfrac{p}{2} z\parcial{}{z}, & p \tx{ is even} \\
x\parcial{}{x}+y\parcial{}{y}+\left(\frac{p+q}{2}\right)
z\parcial{}{z},& p \tx{ and } q \tx{ are odd}.
\end{array}\right.
\]
So, we have that the Hopf fibration $\H_{\F_{\Omega}} \: (f:M\to
\tilde{D})$, adapted to the foliation defined by the one-form
$\Omega\in \sum_{pq}$ will be determined (not uniquely) by a
linearizable singularity of a holomorphic vector  field
$X=X_1+X_2+\cdots$.

Having defined a  Hopf fibration adapted to $\F_{\Omega}$, we can
define the holonomy of the leaf $\tilde{D}\setminus
Sing(\tilde{\F}_{\Omega})$ respect to this fibration. In order to
determine it, we fix a point $p_0\in \tilde{D}\setminus
Sing(\tilde{\F}_{\Omega})$. Over this point we have a transversal
$f^{-1}(p_0)$ and by  path lifting construction, a representation of
the fundamental group of $\tilde{D}\setminus
Sing(\tilde{\F}_{\Omega})$ in $\diff (\C,0)$ is determined, denoted
by ${\mathcal H}_{\Omega,\tilde{D}}$
\[
{\mathcal
H}_{\Omega,\tilde{D}}:\pi_1(\tilde{D}\setminus
Sing(\tilde{\F}_{\Omega}),p_0)\to Diff(\C,0).
\]
This representation is independent of $p_0$ modulo conjugacy and its
image will
be called the \emph{exceptional holonomy} and denoted
$H_{\Omega,\tilde{D}}$.

\section{Classification of the singularities}

From  section (\ref{section3}) we know that the homotopy group
$\pi_1(\tilde{D}\setminus Sing(\tilde{\F}_{\Omega});t_0)$ can be
generated by two elements $\alpha$ and $\beta$ in all the cases
considered, with different relations in each case:
\begin{enumerate}
\item If $p$ and $q$ are even:
$\alpha^{\frac{q}{2}}\beta=\beta\alpha^{\frac{q}{2}}$,

\item If $p$ is even and $q$ is odd: $\alpha^q=\beta^{2}$,

\item If $p$ and $q$ are odd: $\alpha^2\beta=\beta\alpha^2$.
\end{enumerate}

If $\gamma$ is an element of the homotopy group, let us denote
$h_\gamma$ its image by the map ${\mathcal H}_{\Omega, \tilde{D}}$
in the exceptional holonomy. This holonomy can be generated by
$h_{\alpha}$, $h_{\beta}$, which at least satisfy the same relations
than $\alpha$, $\beta$. But in some cases, these relations may be
improved. The following proposition collects some of these
improvements:

\begin{prop}
\begin{enumerate}
\item If $p$ is even,
$h_{\alpha}^{\frac{p}{2}}=id$.
\item If $p$ is even and $q$ is odd,
$h_{\alpha}^{\frac{p}{2}}=h_{\beta}^{p'}=id$.
\end{enumerate}

\end{prop}

\begin{proof}

Consider $p$ even. After $\dfrac{p}{2}$ blow-ups, the
strict transform of the separatrix ${\mathcal S}$ is given by a
surface analytically equivalent to $t^2_{p/2}+y^q=0$. This singular
surface is a cylinder over a curve, that
is either a cuspidal curve of characteristic pair $(2,q)$ or a
couple
of regular curves tangent at the origin at order $\frac{q}{2}$.
Applying Picard-Lefschetz techniques, it can be seen that
the loop $\alpha$ is a simple curve contained in the plane
$y=\varepsilon$, with $|\varepsilon |$ small  enough, that turns
around the points $(t,y)=(\pm i\cdot
 \eps^{q/2}, \eps )$. Thus, the holonomy $h_{\alpha}$ is
completely  determined
by the holonomy of a loop  that turns around the line
$D_{p/2-1}\cap
D_{p/2}$. Along this line, the foliation is a reduced foliation
of dimensional type 2 (in fact, we are in presence of a Kupka
phenomenon), and its analytic type is determined by a
two-dimensional section transversal to the $y$-axis. This
foliation has a linearizable, periodic holonomy, and $h'_{\alpha
} (0)= e^{-2\pi i \cdot \frac{p-2}{p}}$.

If, moreover, $q$ is odd, the periodicity of $h_{\alpha}$ implies
the periodicity of $h_{\beta}$, and so, $h_{\beta}$ is linearizable,
$h'_{\beta}(0)= e^{2\pi i \frac{q'}{p'}}$. Nevertheless, it does not
mean necessary that the holonomy group $H_{\Omega, \tilde{D}}$ is
linearizable, as in particular we don't know if it is abelian or
not.
\end{proof}

%\begin{rem}
%Dar informacion de estos grupos. (Loray).
%\end{rem}

The following theorem contains the main result of the paper. In the
proof, several techniques from \cite{CeMo,BMS,CM,Me} are frequently
used, and we shall not enter in details about them.

\begin{teo} Let $\Omega_1$, $\Omega_2$ be elements of
$\sum_{pq}$. Consider the foliations $\F_{\Omega_1}$ and
$\F_{\Omega_2}$, and their exceptional holonomies
$H_{\Omega_i,\tilde{D}}= <h_{\alpha}^i, h_{\beta}^i >$, $i=1,2$,
defined as before. Then, the foliations are analytically
conjugated if and only if the couples $(h_{\alpha}^i, h_{\beta}^i
)$ are also analytically conjugated, i.e., if and only if there
exists $\Psi \in \diff (\C, 0) $ such that $\Psi^\ast
h^1_{\gamma}= h^2_{\gamma }$, where $\gamma =\alpha, \beta$.
\end{teo}

\begin{prueba} If the foliations are conjugated then clearly
their exceptional holonomies are also conjugated. Conversely,
suppose that the exceptional holonomies are conjugated via $\Psi$,
and let $\tilde{\F}_{\Omega_1}$, $\tilde{\F}_{\Omega_2}$ be the
desingularized, reduced foliations. Because of the existence of
the Hopf fibration relative to $\tilde{D}$, $\Psi$ can be extended
to a neighbourhood of $\tilde{D}$, away from the singular points.
These singular points are the intersections of $\tilde{D}$ with
the other components of the divisor, and with the separatrix (the
separatrices in the even-even case). All these points are singular
points of dimensional types two or three, and for all of them,
$\tilde{D}$ is a strong separatrix. In this situation, the
conjugation of the holonomies of $\tilde{D}$ implies conjugation
of the reduced foliations in a neighbourhood of the singular
points \cite{CeMo}.

So, we have that $\tilde{\F}_{\Omega_1}$, $\tilde{\F}_{\Omega_2}$
are conjugated in a neighbourhood of $\tilde{D}$. Suppose now that
$p$ is even. We need to conjugate the foliations also in a
neighbourhood of $D_1,\ldots D_{p/2-1}$. As $D_1$ is simply
connected, its holonomy is trivial. So, the holonomy of $D_2$,
generated by one loop around $L_1=D_1\cap D_2$ is periodic (the
argument is the same as in \cite{MM}). The same argument shows
that $D_i$ has a periodic holonomy, $1\leq i<\dfrac{p}{2}$, and
so, the foliations have first integrals in a neighbourhood of each
$L_i$, $1\leq i<\dfrac{p}{2}$. These are points of dimensional
type two. By analogous reasons  as in the two-dimensional case,
$\Psi$ can be extended to a neighbourhood of the exceptional
divisor, so, $\F_{\Omega_1}$, $\F_{\Omega_2}$ are conjugated
outside the singular locus, which has codimension two. We conclude
using Hartogs' theorem to extend the conjugation to a
neighbourhood of the origin.

Suppose now that $p$, $q$ are odd. $\tilde{\F}_{\Omega_1}$ and
$\tilde{\F}_{\Omega_2}$ are conjugated in a neighbourhood of
$\tilde{D}$ (that is a projective ${\mathbb P}^2_\C$ in this case).
The fundamental group of $D_{\frac{p-1}{2}}$ is generated by only
one loop, that, after the resolution, can be seen as a loop around
$D_{\frac{p-1}{2}} \cap D''_{(1)}$. This is one of the loops that
generates the holonomy of $D_{\frac{p-1}{2}}$ locally at the reduced
singular points, and following similar arguments as in the preceding
cases, and as the two-dimensional case, the foliation is
linearizable around these points. Let us detail, in this case, how
the use of first integrals allows the extension of the conjugation.

Consider, for instance, the singular point $D_{\frac{p-1}{2}}\cap
D_{\frac{p+q}{2}-2}\cap D_{\frac{p+q}{2}-1 }$, with coordinates
$(x',s',t')$ as in picture \ref{coordinates} .
\begin{figure}[h]
\begin{center}
\includegraphics[height=8cm,width=12cm]{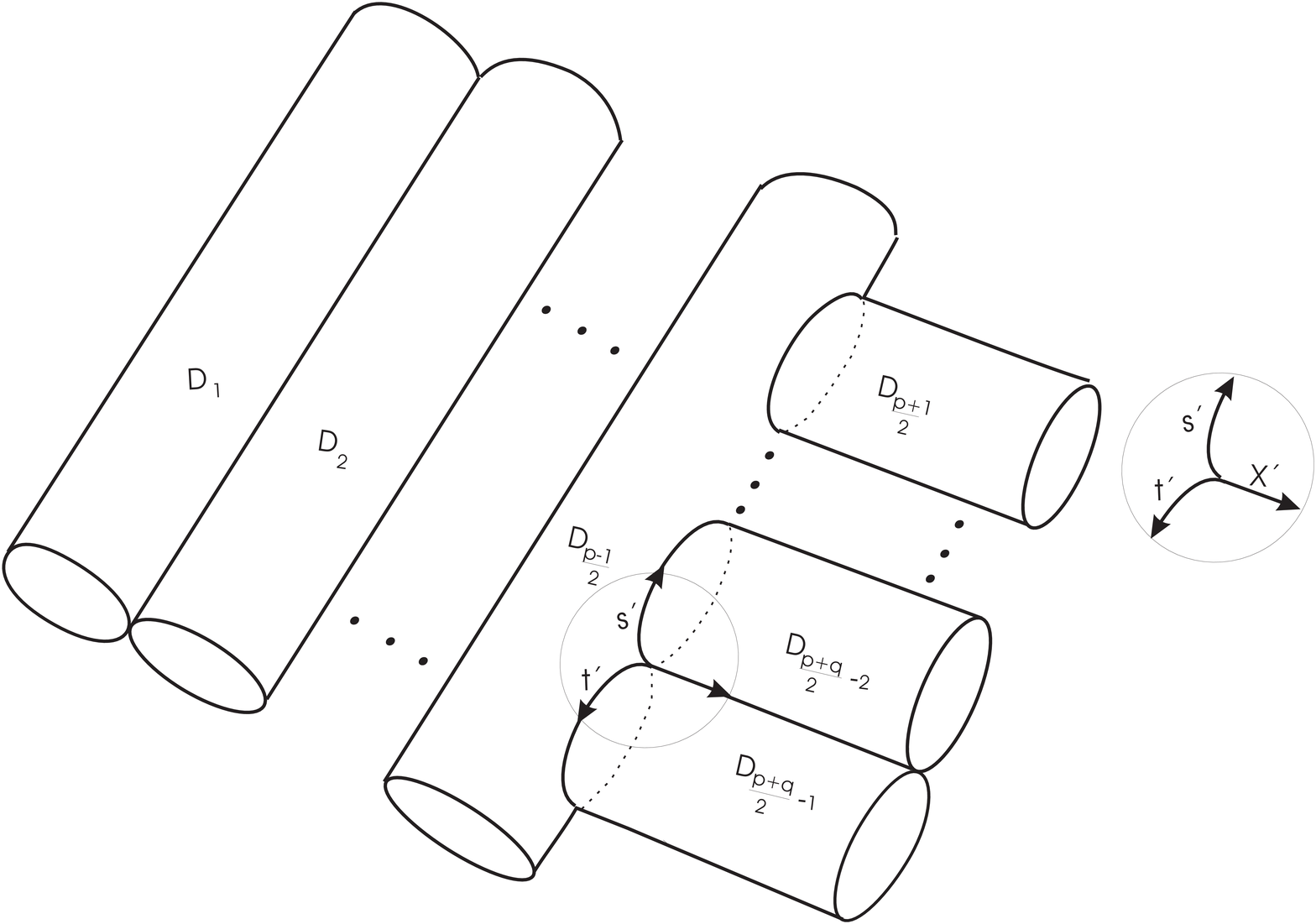}
\caption{\footnotesize Coordinates $(x',s',t')$ the singular point
$D_{\frac{p-1}{2}}\cap D_{\frac{p+q}{2}-2}\cap D_{\frac{p+q}{2}-1
}$}. \label{coordinates}
\end{center}
\end{figure}

We have a conjugation $\Psi$ between the foliations
$\F_{\Omega_1}$ and $\F_{\Omega_2}$ defined over an annulus
$$
\{ |x'| <\eps \} \times \{ |s'| <\eps \} \times \{ c_1<|t'|
<c_2 \}
$$
that respects the fibration. In these coordinates, the foliation
is given by $x't'=cst.$, $s't'^2=cst.$, and the first integral of
$\F_{\Omega_j}$ is $x'^{p-1}s'^{q-1}t'^{q-3} \cdot U_j
(x',s',t')$, $U_j(0)=1$.  This first integral may be extended to
$$
\{ |x'| <\eps \} \times \{ |s'| <\eps \} \times \{ |t'|
<c \} ,
$$
where $c_1<c$, eventually making $\eps$ small enough. We look
first for a diffeomorphism $\Psi_j$ that transforms this first
integral into $x'^{p-1}s'^{q-1}t'^{q-3}$, respecting the
fibration. This diffeomorphism is
$$
\Psi_{j}(x',s',t')= (x'\cdot V_{1j}, s'\cdot V_{2j}, t'\cdot
V_{3j}),
$$
and the conditions mean that
\begin{eqnarray*}
V_{1}^{p-1}\cdot V_2^{q-1}\cdot V_3^{q-3} & = & U_j \\
V_{1}\cdot V_3 & = & 1 \\
V_2\cdot V_3^2 & = & 1.
\end{eqnarray*}
So, $V_{1j}=U_j^{-(p+q)}$; $V_{2j}= U_j^{-2(p+q)}$;
$V_{3j}=U_j^{p+q}$.

Consider now the diffeomorphism $\tilde{\Psi}:= \Psi_{1}\circ
\Psi\circ \Psi_2^{-1}$. It respects both the fibration and the
first integral $x'^{p-1}s'^{q-1}t'^{q-3}$. Write $\tilde{\Psi}=
(\theta_1, \theta_2, \theta_3 )$. The conditions above mean that
\begin{eqnarray*}
\theta_1\cdot \theta_3 & = & x't' \\
\theta_2 \cdot \theta_3^2 & = & s't'^2 \\
\theta_1^{p-1} \cdot \theta_2^{q-1} \cdot \theta_3^{q-3} & = &
x'^{p-1}s'^{q-1}t'^{q-3} \cdot g(x'^{p-1}s'^{q-1}t'^{q-3}),
\end{eqnarray*}
with $g(0)\neq 0$. As before, we have that $\theta_1 =x'\cdot
g^{-(p+q)}$; $\theta_2 = s' \cdot g^{-2(p+q)}$; $\theta_3 =
t'\cdot g^{p+q}$.

This is a map defined, in the considered chart, over a set of the
type $\{ |x'^{p-1}s'^{q-1}t'^{q-3} | <\eps \}$, and this set
intersects the domain of definition of $\Psi$. So, $\Psi=
\Psi_1^{-1}\circ \tilde{\Psi}\circ \Psi_2$ may be extended to a
neighbourhood of $L_{\frac{p+q}{2}-1}= D_{\frac{p-1}{2}}\cap
D_{\frac{p+q}{2}-1}$.

Repeating the argument, we extend the conjugation to a neighbourhood
of $D_{\frac{p-1}{2}}\cap (D_{\frac{p+1}{2}}\cup \cdots \cup
D_{\frac{p+q}{2}})$. Now, similar arguments as in the preceding
si\-tua\-tions, and as in the two-dimensional case, allow us to
extend $\Psi$ to a neighbourhood of the exceptional divisor, and
again Hartogs' theorem completes the result.

\end{prueba}

%-------------------------------------------------------------------------------------------

% Referencias

%-------------------------------------------------------------------------------------------

\end{document}